\documentclass{amsart}
\usepackage[top=1.15in, bottom=1.15in, left=1.17in, right=1.17in]{geometry}
\usepackage{amssymb}
\usepackage{pdflscape}
\usepackage{amscd}
\usepackage{fancyhdr}
\usepackage{tikz-cd}

\newcommand{\naturals}{\mathbb{N}}
\newcommand{\integer}{\mathbb{Z}}
\newcommand{\rational}{\mathbb{Q}}
\newcommand{\divides}{|}
\newcommand{\complex}{\mathbb{C}}
\newcommand{\homs}{{\rm Hom}}
\newcommand{\acts}{\cdot}
\newcommand{\tr}{{\rm tr}}

\newtheorem{theorem}{Theorem}[section]

\newtheorem{lemma}[theorem]{Lemma}
\newtheorem{proposition}[theorem]{Proposition}

\newtheorem{example}[theorem]{Example}
\newtheorem{definition}[theorem]{Definition}

\begin{document}



\title{Donaldson--Thomas invariants of symmetric quivers}

\author{Markus Reineke}
\address{Faculty of Mathematics, Ruhr University Bochum, D-44780 Bochum, Germany}

\maketitle

\begin{abstract}
We review several algebraic, combinatorial and geometric interpretations of motivic Donaldson--Thomas invariants of symmetric quivers.
\end{abstract}

\section{Introduction} 
Motivic Donaldson--Thomas (DT) invariants are introduced in \cite{KSDT,KSCoha} as a mathematical definition of BPS state counts. Though one is mainly interested in these invariants for $3$--Calabi--Yau varieties, or their local versions of $3$--Calabi--Yau categories of representations of quivers with potential and stability, the simplest case of DT invariants of symmetric quivers already produces a rich theory.

The latter have the advantage of admitting an elementary definition, via a product factorization of series of $q$-hypergeometric type. These invariants were proven to be polynomials with nonnegative coefficients in \cite{E}. It is thus natural to ask for interpretations of these coefficients as dimensions of algebraic objects, in particular as dimensions of cohomology groups of geometric objects, or as counts of combinatorial objects, associated to the quiver. In this note, we give a short survey of several such interpretations. 

We start in Section \ref{s2} with the elementary definition and compile a list of examples with which the interpretations of the following sections will be illustrated. Representations of quivers are briefly reviewed in Section \ref{s3}, and the motivic generating series is introduced. Section \ref{s4} surveys geometric interpretations of DT invariants. We introduce moduli spaces of (semi-)simple representations of symmetric quivers, and interpret DT invariants in terms of their intersection cohomology and Chow groups, respectively, following \cite{FR,MeR}. We also describe the interpretation in terms of cohomology of Nakajima quiver varieties of \cite{HLV}. 

Turning to algebraic interpretations of DT invariants in Section \ref{s5}, we first define the Cohomological Hall algebra (CoHa) of a quiver following \cite{KSCoha}. We formulate Efimov's theorem \cite{E} on the structure of the CoHa. We then discuss two dual pictures: first, the graded dual of the CoHa is the enveloping vertex algebra of a vertex Lie algebra \cite{DM}, second, the Koszul dual of an explicit algebra introduced in \cite{DFR} is the enveloping superalgebra of a Lie superalgebra. Both Lie algebras again encode the DT invariants.

Finally, we describe two combinatorial formulas for DT invariants (more precisely, their specializations at $1$). The first identifies the numerical DT invariant as the count of so-called break divisors on a covering graph of the symmetric quiver following \cite{RRT}. The second formula essentially reduces calculation of the DT invariants to those of multiple loop quivers of \cite{RdCoha}, by inductively ``diagonalizing'' the quiver while keeping track of the behaviour of the invariants.\\[2ex]
{\bf Acknowledgments:} The material of this survey is based on lectures at the workshops ``DT--invariants of symmetric quivers: algebra, combinatorics and topology'' (University of Strasbourg) and ``Learning workshop: Knots, homologies and physics'' (University of Warsaw). The author would like to thank the organizers of both workshops for the inspiring atmosphere, and Lydia Gösmann for careful reading of a preliminary version.

\section{Elementary definition of Donaldson--Thomas invariants of symmetric quivers}\label{s2}

The motivic Donaldson--Thomas invariants of symmetric quivers admit a definition without reference to quantum field theory, algebraic geometry or representation theory, by factoring a $q$-hypergeometric series into an infinite product.

To a symmetric matrix, we associate a $q$-hypergeometric function as follows \cite{KSCoha}:

\begin{definition}\label{DTA} Let \(A \in M_{n \times n}(\naturals)\) be a symmetric matrix with non-negative integer entries. Define a formal series \(P_A(x)\in\mathbb{Q}(q^{1/2})[[x_1,\ldots,x_n]]\) of $q$-hypergeometric type as follows:

\[P_A(x) = \sum_{d_1, \ldots, d_n \geq 0} \frac{(-q^{{1}/{2}})^{\sum_{i, j} (a_{ij} - \delta_{ij}) d_i d_j}}{\prod_{i=1}^n(1-q^{-1})\cdot\ldots\cdot(1-q^{-d_i})} x_1^{d_1}\cdot\ldots\cdot x_n^{d_n}.\]
\end{definition}

In the following, we will abbreviate $x_1^{d_1}\cdot\ldots\cdot x_n^{d_n}$ to $x^{\bf d}$ and always understand ${\bf d}$ to be a multiindex of nonnegative integers.

Note that any formal series in $\mathbb{Q}[[x]]$ with constant term \(1\) can be written as an infinite product
\[1 + \sum_{n\geq 1}c_nx^n = \prod_{n \geq 1} (1-x^n)^{-d_n}\]
(expand the right hand side using geometric series and compare coefficients).

Generalizing this remark to the above ring of formal series, we can define:

\begin{definition} Write  \(P_A(x)\) as
\[P_A(x) = \prod_{{\bf d} \neq 0} \prod_{i \in \integer} \prod_{k \geq 0} (1-q^{k+\frac{i+1}{2}}x^{\bf d})^{(-1)^{i+1}c_{{\bf d}, i}}.\]
The {Donaldson--Thomas (DT) invariants} of \(A\) are then defined as the formal Laurent series
\[DT_{\bf d}^A(q) = \sum_{k \in \integer} c_{{\bf d}, k} (-q^{{1}/{2}})^k \in \rational[[q^{\pm{1}/{2}}]]\]
for ${\bf d}\not=0$.  We also define normalized Donaldson--Thomas invariants by
$$\widetilde{DT}_{\bf d}^A(q)=(-q^{1/2})^{\sum_{i, j} (a_{ij} - \delta_{ij}) d_i d_j + 1} DT_{\bf d}^A(q).$$
\end{definition}

The additional product over the index $k$ in the factorization of $P_A(x)$ will become natural by the following examples. The choices of signs in motivic Donaldson--Thomas theory are always subtle and an essential feature. The main justification for these choices is the following theorem confirming a conjecture of \cite{KSCoha}:

\begin{theorem}[\cite{E}]\label{E} For all $A$ and ${\bf d}$ as above, the normalized DT invariant is a polynomial in $-q^{1/2}$ with nonnegative coefficients:
 $$\widetilde{DT}_{\bf d}^A(q) \in \naturals[-q^{{1}/{2}}].$$
\end{theorem}

In the light of this theorem, it is natural to ask for interpretations of the numbers $c_{{\bf d},k}$ as dimensions of vector spaces of algebraic objects associated to $A$, in particular dimensions of cohomology groups, or as counts of combinatorial objects. Such objects will be associated to (representations of) quivers, using an interpretation of DT invariants in terms of the motivic generating series of a quiver.

\begin{example}\label{ex1} We first work out three cases where all DT invariants for a matrix $A$ can be calculated. We then proceed by stating several more explicit examples from the literature, which will be re-derived using the various interpretations of DT invariants to be described in the following sections.

 \begin{enumerate}
 
 \item For the first example, we use the $q$-binomial theorem:
$$\sum_{k\geq 0}\frac{t^k}{(1-q)\cdot\ldots\cdot(1-q^k)}=
\prod_{k\geq 0}(1-q^kt)^{-1}$$ 

  For the matrix \(A=[0]\), we have
  \[P_A(x) = \sum_{d \geq 0} \frac{(-q^{{1}/{2}})^{-d^2}}{(1-q^{-1})\cdot\ldots\cdot(1-q^{-d})} x^d.\]
 Expanding the fraction by \(q^{d(d+1)}/{2}\) and applying  the $q$-binomial theorem, this reads
  \[P_A(x) = \sum_{d \geq 0} \frac{(q^{{1}/{2}}x)^d}{(1-q) \cdot\ldots\cdot (1-q^d)}=\prod_{k \geq 0} (1-q^{k+{1}/{2}}x)^{-1}.\]
  Thus we find \[DT_1^{A}(q) = 1,\; DT_{\geq 2}^{A}(q) = 0.\]
  
\item  Similarly, we treat  the matrix $A=[1]$. We use another version of the $q$-binomial theorem:

$$ \sum_{k\geq 0}\frac{q^{k(k-1)/2}t^k}{(1-q)\cdot\ldots\cdot(1-q^k)}=\prod_{k\geq 0}(1+q^kt).$$

We have
  $$P_A(x)= \sum_{d \geq 0} \frac{x^d}{(1-q^{-1})\cdot\ldots\cdot(1-q^{-d})}.$$
  Again expanding the fraction by $q^{d(d+1)/2}$ and applying the  $q$-binomial theorem, we find
  $$P_A(x)=\sum_{d\geq 0}\frac{q^{d(d-1)/2}(-qx)^d}{(1-q)\cdot\ldots\cdot(1-q^d)}=\prod_{k\geq 0}(1-q^{k+1}x),$$
  and thus \[DT_1^{A}(q) = -q^{1/2},\; DT_{\geq 2}^{A}(q) = 0.\]

\item For the matrix $A=\left[\begin{array}{cc}0&1\\ 1&0\end{array}\right]$, one finds \cite{FR} $$DT^A_{(1,0)}(q)=1=DT^A_{(0,1)}(q),\, DT^A_{(1,1)}(q)=-q^{1/2}$$ and $DT^A_{\bf d}(q)=0$ for all other ${\bf d}$. In fact, these three symmetric nonnegative integer matrices are the only indecomposable ones (that is, they cannot be transformed to proper block diagonal form by a simultaneous permutation of rows and columns) for which there are only finitely many non-vanishing DT invariants, as will become apparent in the following sections.

\item For a more general matrix $A=\left[\begin{array}{cc}0&m\\ m&0\end{array}\right]$ with $m\geq 1$ and ${\bf d}=(1,k)$, one finds \cite{FR} $$DT^A_{\bf d}(q)=(-q^{1/2})^{k^2}\frac{(1-q)\cdot\ldots\cdot(1-q^m)}{(1-q)\cdot\ldots\cdot (1-q^k)\cdot(1-q)\cdot\ldots\cdot(1-q^{m-k})}$$ for $k\leq m$ and $DT^A_{\bf d}(q)=0$ for $k>m$.

\item As a final example \cite{RdCoha}, for $A=[m]$ with \(m \geq 2\), we have
\[  {\widetilde{DT}_d^{A}}(1) = \frac{1}{d^{2}}\sum_{e \divides d} \mu(\frac{d}{e}) (-1)^{(m-1)(d-e)} \binom{me-1}{e-1}.\]

Note that this is indeed a non-negative integer, but for nontrivial reasons.

\item Trying to extend the definition of $P_A(x)$ to not necessarily symmetric matrix, we see that this function only depends on the values $a_{ij}+a_{ji}$. In case these sums are odd, we cannot expect good properties of the DT invariants, as the following example shows. So there is indeed no loss of generality in assuming $A$ to be symmetric. 
Indeed, for $A=\left[\begin{array}{cc}0&1\\ 0&0\end{array}\right]$ the first nontrivial example $DT_{(1,1)}(q)=-(q+q^{1/2})/(1-q)$ already is a formal series in $-q^{1/2}$ with both positive and negative coefficients.

\end{enumerate}

\end{example}

\section{Quiver representations and their motivic generating functions}\label{s3}

We introduce basic concepts of the representation theory of quivers \cite{S}:

\begin{definition}\label{quiver} A {quiver} $Q$ is an oriented graph, encoded in a set of vertices $Q_0$, and a set of arrows denoted by $\alpha:i\rightarrow j$. We call $Q$ symmetric if the number of arrows from $i$ to $j$ equals the number of arrows from $j$ to $i$ for all $i,j\in Q_0$.  To a matrix \(A\) with nonnegative integer coefficients as before, we associate a symmetric quiver \(Q^A\) with $Q^A_0=\{1,\ldots,n\}$ and number of arrows from $i$ to $j$ equal to $a_{ij}$, for all $i,j\in Q^A_0$.

A representation $V$ of $Q$ consists of  a finite-dimensional complex vector space \(V_i\) for any $i\in Q_0$, and a \(\complex\)-linear map \(f_\alpha:V_i \rightarrow V_j\) for any arrow \(\alpha:i \to j\). The tuple ${\bf d}=(\dim V_i)_{i\in Q_0}$ is called the dimension type of $V$.

Two representations \(((V_i), (f_\alpha))\) and \(((V_i'), (f_\alpha'))\) are equivalent if there are isomorphisms \(\varphi_i \colon V_i \to V_i'\) inducing commutative squares at all \(\alpha\):
\[
  \begin{tikzcd}
    V_i \arrow{r}{f_\alpha} \arrow{d}{\varphi_i}
    & V_j \arrow{d}{\varphi_j} \\
    V_i' \arrow{r}{f_\alpha'}
    & V_j'.
  \end{tikzcd}
\]
A subrepresentation of $V$ consists of a tuple of subspaces $(U_i\subset V_i)_i$ compatible with the structure maps, that is, $f_\alpha(U_i)\subset U_j$ for all $\alpha:i\rightarrow j$. In this case, the quotient spaces and the induced maps $((V_i/U_i)_i,(\overline{f_\alpha})_\alpha)$ form $V/U$, the quotient representation of $V$ by $U$.

The representation $V$ is called simple if it does not admit proper non-zero subrepresentations. It is called semisimple if it is equivalent to a direct sum of simple ones, where direct sums of quiver representations are formed in the obvious way.
\end{definition}

A complete classification and parametrization of all (finite-dimensional) representations of a quiver $Q$ is known if and only if all connected components of the unoriented graph underlying $Q$ is of (simply-laced) Dynkin type or extended Dynkin type.

The only symmetric and connected such quivers \(Q^A\) are associated to $[0]$ (Dynkin type $A_1$),  $[1]$ (type $\widetilde{A}_0$) and  $\begin{bmatrix} 0 & 1 \\ 1 & 0 \end{bmatrix}$ (type $\widetilde{A}_1$).

Before interpreting (the equivalence of) quiver representations geometrically, we define motives of complex algebraic varieties (see, for example, \cite{B}).

\begin{definition}Let $K_0({\rm Var}_\mathbb{C})$ be the free abelian group in isomorphism classes $[X]$ of (quasi-projective) complex algebraic varieties, modulo the relations $[X]=[U]+[A]$ whenever $U\subset X$ is Zariski-open with complement $A=X\setminus U$. It becomes a ring via $[X]\cdot[Y]=[X\times Y]$. We denote by $\mathbb{L}=[\mathbb{A}^1]$ the class of the affine line, and we consider the localization $$R^{\rm mot}=K_0({\rm Var}_\mathbb{C})[\mathbb{L}^{\pm 1/2},(1-\mathbb{L}^i)^{-1}\, :\, i\geq 1].$$ We define  the virtual motive of a connected variety $X$ as $[X]_{\rm vir}=(-\mathbb{L}^{1/2})^{-\dim X}[X]\in R^{\rm mot}$.
\end{definition}

\begin{example} The most relevant  examples of motives for our considerations are
$$[\mathbb{C}^n]=\mathbb{L}^n,\; [\mathbb{P}^n]=(1-\mathbb{L}^{n+1})/(1-\mathbb{L}),\; [{\rm GL}_n(\mathbb{C})]=\mathbb{L}^{n^2}(1-\mathbb{L}^{-1})\cdot\ldots\cdot(1-\mathbb{L}^{-n})$$
and, for the Grassmannian of $k$-dimensional subspaces of an $n$-dimensional vector space,
$$[{\rm Gr}_k^n]=\frac{(1-\mathbb{L})\cdot\ldots\cdot(1-\mathbb{L}^n)}{(1-\mathbb{L})\cdot\ldots\cdot (1-\mathbb{L}^k)\cdot(1-\mathbb{L})\cdot\ldots\cdot(1-\mathbb{L}^{n-k})}.$$
\end{example}

Representations of a quiver \(Q\) on fixed vector spaces $(V_i)_i$ can be viewed as points in the complex affine space $$R_{\bf d}(Q)=\bigoplus_{\alpha:i\rightarrow j} \homs(V_i, V_j).$$
The reductive complex algebraic group $G_{\bf d}=\prod_{i\in Q_0}{\rm GL}(V_i)$ acts on $R_{\bf d}(Q)$  by a $Q$-graded analogue of conjugation:
\[(g_i)_i \acts (f_\alpha)_\alpha = (g_j f_\alpha g_i^{-1})_{\alpha:i\rightarrow j}.\]
By the definition of equivalence of quiver representations and by the definition of the series $P_A(x)$ we find:

\begin{lemma} The orbits of this action are the equivalence classes of representations of $Q$ of dimension type ${\bf d}$. We have
$$P_A(x)=\sum_{d\in\mathbb{N}^n}\frac{[R_{\bf d}(Q^A)]_{\rm vir}}{[G_{\bf d}]_{\rm vir}}x^{\bf d}.$$
\end{lemma}
Thus, heuristically, $P_A(x)$ is a partition function ``counting'' quiver representations up to equivalence.

\section{DT invariants as cohomology of moduli spaces}\label{s4}

The orbit space $R_{\bf d}(Q)/G_{\bf d}$ a priori has a very bad topology. Geometric Invariant Theory \cite{Mu} is a method to construct such quotients as algebraic varieties, we apply it to our setup (see also \cite{Ki}).

The action of $G_{\bf d}$ on $R_{\bf d}(Q)$ induces an action on the ring $\mathbb{C}[R_{\bf d}(Q)]$ of polynomial functions on $R_{\bf d}(Q)$. By a classical theorem of Hilbert, the subring $\mathbb{C}[R_{\bf d}(Q)]^{G_{\bf d}}$ of $G_{\bf d}$-invariant polynomial functions is finitely generated, and thus there exists a complex affine variety $M_{\bf d}(Q)={\rm Spec}(\mathbb{C}[R_{\bf d}(Q)]^{G_{\bf d}})$ having it as its ring of polynomial functions. 

\begin{theorem} The points of the variety $M_{\bf d}(Q^A)$  parametrize equivalence classes of semisimple representations of $Q^A$ of dimension type ${\bf d}$. The variety $M_{\bf d}(Q^A)$ admits a Zariski-open subset $M_{\bf d}^{\rm simp}(Q^A)$ parametrizing the equivalence classes of simple representations, which is irreducible of dimension $\sum_{i,j}(a_{ij}-\delta_{ij})d_id_j+1$.
The invariant ring $\mathbb{C}[R_{\bf d}(Q^A)]^{G_{\bf d}}$ is generated by traces $t_\omega$ along oriented cycles in $Q^A$.
\end{theorem}

The invariant polynomial functions $t_\omega$ are defined as follows: for any oriented cycle
$$\omega:i_1\stackrel{\alpha_1}{\rightarrow}i_2\stackrel{\alpha_2}{\rightarrow}\ldots\stackrel{\alpha_s}{\rightarrow}i_{s+1}=i_1$$
in $Q$, the function $t_\omega$ maps a representation given by linear maps $(f_\alpha)_\alpha$ to ${\rm Tr}(f_{\alpha_s}\circ\ldots\circ f_{\alpha_1})$.

 The last part of the theorem provides explicit coordinates for the moduli spaces $M_{\bf d}^A(Q)$.

\begin{example} We describe some small moduli spaces explicitly:
\begin{enumerate}
\item We consider the matrix  $A=\left[\begin{array}{cc}0&2\\ 2&0\end{array}\right]$ and the dimension type \({\bf d} = (1, 1)\), thus a representation is given by scalars $(a,b,c,d)$ representing the arrows, on which $\mathbb{C}^*\times\mathbb{C}^*$ acts via
  $$(s,t)\cdot(a,b,c,d)=(\frac{t}{s}a,\frac{t}{s}b,\frac{s}{t}c,\frac{s}{t}d).$$
  
  The ring of invariants is generated by the traces (of $1\times 1$-matrices) along oriented cycles
  \[x_1 = ac, \quad x_2 = ad, \quad x_3 = bc, \quad x_4 = bd,\]
  fulfilling the relation \(x_1 x_4 = x_2 x_3\). We thus find that the moduli space is the conifold singularity, $$M_{\bf d}(Q^A)\simeq\{(x_1,x_2,x_3,x_4)\in\mathbb{C}^4\, :\, x_1x_4=x_2x_3\},$$ and $M_{\bf d}^{\rm simp}(Q^A)$ is the complement of the origin.

\item   For \(A = [2]\), and \(d = (2)\), we denote the $2\times 2$-matrices representing the two loops by \(X\) and \(Y\). Then
  \[\tr(X), \tr(X^2), \tr(Y), \tr(Y^2), \tr(XY)\]
  are known to generate the invariant ring, and to be algebraically independent, and so 
  \(M_2(Q^A) \simeq \complex^5.\) The open subset $M_2^{\rm simp}(Q^A)$ is defined by a certain cubic polynomial in these five traces to be non-zero, expressing that $X$ and $Y$ are not simultaneously triangularizable.
  \end{enumerate}
  \end{example}
  
  The main result is that the DT invariants $DT_{\bf d}^A(q)$ can be interpreted as topological invariants of the spaces $M_{\bf d}(Q^A)$, namely as their Poincar\'e polynomials in intersection cohomology with compact support:

\begin{theorem}[\cite{MeR}]\label{MeR} We have
\[  \widetilde{DT}_{\bf d}^A(q) = \sum_k \dim {\rm IH}^k_c(M_{\bf d}(Q^A))(-q^{{1}/{2}})^k\]
  if $M_{\bf d}^{\rm simp}(Q^A)\not=\emptyset$, and equal to zero otherwise.
\end{theorem}
  
Intersection cohomology is usually extremely difficult to compute. In fact, the above theorem often works in the other direction: the DT invariant can be computed with algebraic or geometric methods, and the theorem allows to compute intersection cohomology from this.

One notable exception is formed by varieties $X$ admitting a small resolution of singularities, namely a proper and birational map $f:Y\rightarrow X$ from a smooth variety $Y$, such that $${\rm codim}_X\{x\in X\,:\,\dim f^{-1}(x)\geq r\}>2r$$ for all $r\geq 1$. In this case, we can compute intersection cohomology of $X$ by the cohomology of $Y$: $${\rm IH}_c^*(X)\simeq H_c^*(Y).$$

\begin{example}\label{ex4} We describe two explicit calculations of DT invariants using Theorem \ref{MeR}, as well as a criterion for their non-vanishing.
\begin{enumerate}
\item For $A=\left[\begin{array}{cc}0&m\\ m&0\end{array}\right]$ and ${\bf d}=(1,k)$ for $k\leq m$, a representation consists of $m$ vectors $v_1,\ldots,v_m$ in $\mathbb{C}^k$ and $m$ covectors $\varphi_1,\ldots,\varphi_m$ on $\mathbb{C}^k$. The $m\times m$-matrix $(\varphi_i(v_j))_{i,j}$ is of rank at most $k$ (and it is of rank precisely $k$ if and only if the vectors and covectors define a simple representation of $Q^A$), so $M_{\bf d}(Q^A)$ is isomorphic to the space of $m\times m$-matrices of rank at most $k$. Its singularities can be resolved by the variety $Y$ of pairs $(B,U)\in M_{\bf d}(Q^A)\times{\rm Gr}_k^m$ such that $U$ contains the image of $B$. Then $Y$ is a homogeneous bundle of rank $km$ over the Grassmannian ${\rm Gr}_k^m$. Consequently, we retrieve the description of the DT invariant from Example \ref{ex1} (see also \cite{FR}).

\item We sketch a generalization of the previous example. Whenever the dimension vector ${\bf d}$ is indivisible, one can define moduli spaces $M_{\bf d}^{\Theta-{\rm st}}(Q^A)$ of so-called $\Theta$-stable representations of $Q^A$ of dimension vector ${\bf d}$ \cite{Ki}. For a generic choice of the so-called stability parameter $\Theta$, these moduli spaces provide small desingularizations of $M_{\bf d}(Q^A)$ by \cite{RSmall}. Thus $\widetilde{DT}^A_{\bf d}(q)$ can be computed from the cohomology of $M_{\bf d}^{\Theta-{\rm st}}(Q^A)$ in this case, which can be done using the wall-crossing formula of \cite{RHNS}. 

As an example, consider $A=\left[\begin{array}{cc}1&1\\ 1&1\end{array}\right]$ and ${\bf d}=(2,3)$. The moduli space of stable representations is $13$-dimensional, and with the above method, we find $\widetilde{DT}^A_{\bf d}(q)=q^{13}+q^{12}+2q^{11}$.

\item We can also apply the above theorem to give a numerical criterion for non-vanishing of a DT invariant. We assume without loss of generality that $A$ is indecomposable (so that $Q^A$ is connected), that $d_i\not=0$ for all $i$, and that $A$ is different from the matrices $[0]$, $[1]$, $\begin{bmatrix} 0 & 1 \\ 1 & 0 \end{bmatrix}$ already treated completely. Then, by \cite{LBP} we have $M_{\bf d}^{\rm simp}(Q^A)\not=\emptyset$ (and, by Theorem \ref{MeR}, $DT^A_{\bf d}(q)\not=0$), if and only if $d_i\leq\sum_ja_{ij}d_j$ for all $i\in Q^A_0$.

\end{enumerate}

\end{example}

An interpretation of DT invariants in terms of another homology-type theory for algebraic varieties, namely (rational) Chow groups $A_*(X)_\mathbb{Q}$, is developed in \cite{FR}.

\begin{theorem}[\cite{FR}] We have
\[  \widetilde{DT}_{\bf d}^A(q) = \sum_{i\in \mathbb{N}} \dim A_i(M_{\bf d}^{\rm simp}(Q^A))_\mathbb{Q}q^i.\]
 \end{theorem}
 
\begin{example} We continue the first of the previous examples. The variety of $m\times m$-matrices of rank $k$ fibres over ${\rm Gr}_k^m$ by associating to a matrix its image, with fibre being the space of $k\times m$-matrices of highest rank. This allows to compute the Chow groups as those of the Grassmannian shifted by $k\cdot m$. We recover the above formula.
 \end{example}
 
A third geometric interpretation is provided by the work of Hausel, Letellier and Rodriguez--Villegas \cite{HLV}. For simplicity we assume that all $a_{ii}$ are even. Then $Q^A$ is the double of a (non-symmetric) quiver, and $R_{\bf d}(Q^A)$ admits a natural algebraic symplectic structure making the $G_{\bf d}$-action Hamiltonian, with associated moment map $$\mu_{\bf d}:R_{\bf d}(Q^A)\rightarrow{\rm Lie}(G_{\bf d})^*.$$ For a sufficiently general coadjoint orbit $\mathcal{O}$ of $G_{\bf d}$ in ${\rm Lie}(G_{\bf d})^*$, there exists a geometric quotient $\mu_{\bf d}^{-1}(\mathcal{O})/G_{\bf d}$ as an algebraic variety, which carries an action of the Weyl group $$W_{\bf d}\simeq \prod_{i\in Q^A_0}S_{d_i}$$ of $G_{\bf d}$. DT invariants are then encoded in the $W_{\bf d}$-invariant part of cohomology of these quotients.

\begin{theorem}[\cite{HLV}] We have 
\[ \widetilde{DT}_{\bf d}^A(q) =q^{\sum_{i, j} (a_{ij} - \delta_{ij}) d_i d_j + 1} \sum_k \dim H^k(\mu_{\bf d}^{-1}(\mathcal{O})/G_{\bf d})^{W_{\bf d}}(-q^{1/2})^{-k}.\]
 \end{theorem}

\section{Cohomological Hall algebras, or algebras of BPS states}\label{s5}

Instead of relating DT invariants to the cohomology of quotient varieties, we can consider the equivariant cohomology of the representation spaces themselves. This leads to the notion of Cohomological Hall algebras of \cite{KSCoha}, as a mathematical definition of the algebra of BPS states of \cite{HM} (see also \cite{FRK} for a short introduction).

\begin{definition} For a quiver $Q$, define its Cohomological Hall algebra (CoHa) as
$$A(Q)=\bigoplus_{{\bf d}}H^*_{G_{\bf d}}(R_{\bf d}(Q)),$$
which carries a natural associative $\mathbb{N}Q_0$-graded multiplication via convolution along a diagram
$$R_{\bf d}\times R_{\bf e}\leftarrow Z_{{\bf d},{\bf e}}\rightarrow R_{{\bf d}+{\bf e}}$$
encoding short exact sequences $0\rightarrow U\rightarrow V\rightarrow V/U\rightarrow 0$ of quiver representations.
\end{definition}

Using localization in equivariant cohomology, one can give an explicit algebraic description of this algebra:

\begin{proposition}[\cite{KSCoha}] We have $$A(Q)=\bigoplus_{\bf d}\mathbb{Q}[x_{i,r}\, :\, i\in Q_0, r=1,\ldots,d_i]^{W_{\bf d}},$$ where the multiplication of $W_{\bf d}$-invariant polynomials is given by a shuffle product with kernel:\\
$(f*g)(x_{i,r}:i\in Q_0, r\leq d_i+e_i)=$
$$=\sum_\sigma f(x_{i,\sigma_i(r)}:i\in Q_0, r\leq d_i)\cdot g(x_{i,\sigma_i(d_i+s)}:i\in Q_0, s\leq e_i)\cdot$$
$$\cdot\prod_{i,j\in Q_0}\prod_{{r\leq d_i}\atop{s\leq e_j}}(x_{j,\sigma_j(d_j+s)}-x_{i,\sigma_i(r)})^{a_{ij}-\delta_{ij}},$$
where the sum ranges over all tuples $\sigma=(\sigma_i)_{i}$ of $(d_i,e_i)$-shuffle permutations.
\end{proposition}

If $Q$ is symmetric, we can use the cohomological grading to refine the grading on  $A(Q)$ to a $\mathbb{N}Q_0\times\mathbb{Z}$-grading, and the product can be modified by a sign to make it graded supercommutative. It is easy to see that the $q$-hypergeometric series $P_A(x)$ then equals the Poincar\'e--Hilbert series of the algebra $A(Q^A)$.

\begin{lemma} We have
$$P_A(x)=\sum_{\bf d}\sum_k\dim A(Q^A)_{{\bf d},k}(-q^{1/2})^{-k}.$$
\end{lemma}

Theorem \ref{E} more precisely reads as follows:

\begin{theorem}[\cite{E}] There exists a $\mathbb{N}Q^A_0\times\mathbb{Z}$-graded subspace $V^{\rm prim}\subset A(Q^A)$, with each graded component $V^{\rm prim}_{{\bf d},*}$ finite-dimensional,  and an element $z\in A(Q^A)$, such that
$$A(Q^A)\simeq{\rm Sym}(V^{\rm prim}[z]).$$
Consequently, we have
$${DT}^A_{\bf d}(q)=\sum_k(-q^{1/2})^{1-k}\dim V_{{\bf d},k}.$$

\end{theorem}

\begin{example} We can treat the cases of the matrices $A=[0]$ and $A=[1]$ explicitly.
\begin{enumerate}
\item For $A=[0]$ the CoHa is a free anticommutative algebra in $\mathbb{Q}[x_{i,1}]$, where the generator $x_{i,1}^k$ is placed in bidegree $(1,2k+1)$. Consequently, we have $V^{\rm prim}=\mathbb{Q}\cdot x_{i,1}$.
\item For $A=[1]$, the CoHa is a free commutative algebra in $\mathbb{Q}[x_{i,1}]$, where the generator $x_{i,1}^k$ is placed in bidegree $(1,2k)$. Consequently, we have $V^{\rm prim}=\mathbb{Q}\cdot 1$.
\end{enumerate}Both descriptions again confirm the values of DT invariants computed in Example \ref{ex1}.
\end{example}

The space $V^{\rm prim}$ of the previous theorem is not canonically defined. So it is desirable to find additional structures on Cohomological Hall algebras making this space canonical. This is not achieved for $A(Q)$ itself, but for various duals which we will now discuss (but we do not give all details about the necessary shifts in the gradings).

First, we consider the graded dual $$A(Q)^\vee=\bigoplus_{\bf d}\bigoplus_kH_{G_{\bf d}}^k(R_{\bf d}(Q))^*.$$

\begin{theorem}[\cite{DM}] The graded coalgebra $A(Q^A)^\vee$ carries a canonical structure of a vertex bialgebra, so that it is isomorphic to the enveloping vertex algebra of its vertex Lie algebra $C$ of primitive elements. The coefficients of the DT invariants are the dimensions of the graded components $C/\partial C$, where $\partial$ denotes the canonical derivation on $C$.
\end{theorem}

Partially building on this vertex algebra picture, we can define an explicit superalgebra by generators and relations, whose Koszul dual turns out to be the enveloping superalgebra of its Lie superalgebra of primitive elements,  again encoding the DT invariants.

\begin{definition}[\cite{DFR}]\label{bqa} Define $B(Q^A)$ as the $\mathbb{N}Q^A_0\times\mathbb{Z}$-graded supercommutative algebra with generators $\alpha_{i,k}$ for $i\in Q^A_0$ and $k\geq 0$, and relations
$$\alpha_i(z)\frac{d^p}{dz^p}\alpha_j(z)=0\mbox{ for all }0\leq p<a_{i,j},$$
where $\alpha_i(z)=\sum_{k\geq 0}\alpha_{i,k}z^k$.
\end{definition}

For an algebra $A$ generated by a vector space $V$ subject to a space of quadratic relation $R\subset V\otimes V$, we define the Koszul dual algebra $A^!$ as being generated by the (graded) dual space $V^\vee$ subject to the relations $R^\perp\subset (V^\vee)\otimes(V^\vee)\simeq (V\otimes V)^\vee$.  For example, the symmetric algebra of a vector space Koszul-dualizes to the exterior algebra of the dual space.

\begin{theorem}[\cite{DFR}]  The Koszul dual of $B(Q^A)$ is the enveloping superalgebra of a graded Lie superalgebra $\mathfrak{g}^A$, carrying an action of the first Weyl algebra $\langle t,\partial_t\rangle$ on each $\mathbb{N}Q^A_0$-graded component. The coefficients of the DT invariants are the dimensions of the spaces ${\rm Ker}(\partial_t)$.
\end{theorem}

In fact, the Lie superalgebra $\mathfrak{g}^A$ also admits an explicit description by generators and relations.

\begin{example} In the case $A=\left[\begin{array}{cc}0&1\\ 1&0\end{array}\right]$, the Lie superalgebra $\mathfrak{g}(Q^A)$ has odd generators $b_{1,k},b_{2,k}$ for $k\geq 0$, subject to the relations $$ [b_{1,k},b_{1,l}]=0,\; [b_{1,k},b_{2,l}]=[b_{1,k-1},b_{2,l+1}],\; [b_{2,k},b_{2,l}]=0.$$
By the second group of relations, the commutators $[b_{1,k},b_{2,l}]$ depend only on $k+l$ and are denoted by $c_{k+l}$. They turn out to be central by the (super) Jacobi identity, therefore $\mathfrak{g}(Q^A)$ has a basis consisting of the $b_{1,k}, b_{2,k}, c_k$ for $k\geq 0$. The derivation $\partial_t$ acts by lowering the index $k$, thus ${\rm Ker}(\partial_t)$ is spanned by $b_{1,0}$, $b_{2,0}$, $c_0$ in degrees $((1,0),1)$, $((0,1),1)$, $((1,1),2)$, respectively, compatible with the DT invariants computed in Example \ref{ex1}.
\end{example}

\section{Combinatorial formulas}\label{s6}

We will state a combinatorial and manifestly positive formula for $\widetilde{DT}^A_{\bf d}(1)$ in graph-theoretic terms, under the assumption that $Q^A$ has at least one loop at any vertex. We then say that $Q^A$ has enough loops; that this assumption is rather natural in the present context can, for example, be seen from the algebraic description of the CoHa: if $Q^A$ has enough loops, there are no proper denominators in the modified shuffle product.

\begin{definition} For an unoriented graph $\Gamma$ with set of vertices $V(\Gamma)$ and set of edges $E(\Gamma)$ we denote by $g(\Gamma)=|E(\Gamma)|-|V(\Gamma)|+1$ its genus. An effective divisor on $\Gamma$ is a function $D:V(\Gamma)\rightarrow\mathbb{N}$. Its degree $d(D)$ is defined as $\sum_{v\in V(\Gamma)}D(v)$; the divisor $D$ can be restricted to a divisor $D|_{\Gamma'}$ on any subgraph $\Gamma'\subset\Gamma$. We call $D$ a break divisor if $d(D)=g(\Gamma)$, and $d(D|_{\Gamma'})\geq g(\Gamma')$ for any connected $\Gamma'\subset\Gamma$. 
\end{definition}

\begin{theorem}[\cite{RRT}] Given $Q^A$ with enough loops and ${\bf d}\in\mathbb{N}Q^A_0$, define $\Gamma_{\bf d}(Q^A)$ as the graph with vertices $v_{i,k}$ for $i\in Q^A_0$ and $1\leq k\leq d_i$, and $a_{i,j}-\delta_{i,j}$ edges between $v_{i,k}$ and $v_{j,l}$ for all $(i,k)\not=(j,l)$. The group $W_{\bf d}$ acts naturally on $\Gamma_{\bf d}(Q^A)$ by permuting vertices associated to the same $i$. Then $\widetilde{DT}^A_{\bf d}(1)$ equals the number of $W_{\bf d}$-orbits of break divisors on $\Gamma_{\bf d}(Q^A)$.
\end{theorem}

As an example, consider $A=\left[\begin{array}{cc}1&1\\ 1&1\end{array}\right]$ and ${\bf d}=(2,3)$. Then $\Gamma_{\bf d}(Q^A)$ is the complete bipartite graph with groups of vertices $v_{1,1}, v_{1,2}$ and $v_{2,1}, v_{2,2}, v_{2,3}$, respectively. A break divisor is thus given by a tuple $(a_1,a_2,b_1,b_2,b_3)$ such that $a_1+a_2+b_1+b_2+b_3=2$ and $a_1+a_2+b_k+b_l\geq 1$ for all $1\leq k<l\leq 3$. These are considered up to reindexing of the $a$- and $b$-variables. We arrive at four representatives $$(2,0,0,0,0), (1,1,0,0,0), (1,0,1,0,0), (0,0,1,1,0),$$ thus $\widetilde{DT}^A_{\bf d}(1)=4$, compatible with Example \ref{ex4}.

Finally, we discuss quiver diagonalization \cite{diag}. For this, we extend the definition of the series $P_A(x)$ to infinite matrices over an index set $J$.

\begin{theorem} For every symmetric matrix $A$ there exists a (infinite) diagonal matrix $A^\infty$, together with monomials $M_j(q,x_1,\ldots,x_n)$, such that $P_A(x)$ equals $P_{A^\infty}(x^\infty)$ under the specialization of variables $x^\infty_j\mapsto M_j(q,x)$. 
\end{theorem}

Consequently, this allows to write DT invariants $DT^A_{\bf d}(q)$ as a sum of invariants $DT^{A^\infty}_{{\bf d}'}(q)$ for various dimension vectors ${\bf d'}$ for $Q^{A^\infty}$. Since $A^\infty$ is diagonal, the latter invariants can be computed with any of the techniques developed so far, for example their value at one with the formula from Example \ref{ex1}. For more details, we refer to \cite{No}. An interpretation of the quiver diagonalization process in terms of the algebra $B(Q^A)$ above is given in \cite{DFKR}.

\end{document}